\documentclass{article} 
\usepackage{amsmath}
\usepackage{amssymb}
\usepackage{amsxtra}
\usepackage{amsthm}
\numberwithin{equation}{section}
\newcommand{\comment}[1]{}

\newcommand\hs[3]{\left(
\begin{array}{c}{#1}\\[6pt]{#2}\end{array}\,;\,{#3}\,\right)
}
\allowdisplaybreaks[1]

\newcommand\C{{ \Bbb C}}

\newcommand\PP{{\Bbb P}}

\newcommand\e{{\varepsilon}}

\newtheorem{lem}{Lemma}

\newtheorem{prop}{Proposition}

\begin{document}
\begin{flushleft}
{\Large\bf Regularizable cycles 
associated with a Selberg type integral 
under some resonance condition}
\bigskip
\bigskip

{Katsuhisa Mimachi and Masaaki Yoshida}
\end{flushleft}

\bigskip

{ \bf Abstract.}  We study the twisted homology group
attached to a Selberg type integral under some resonance 
condition, which naturally appears in the $su_2$-conformal
field theory and the representation of the Iwahori-Hecke
algebra. We determine the dimension of the
space of the regularizable cycles. The dimension-formula is given
in terms of the generalized hypergeometric series ${}_3F_2.$
\vspace*{1cm}

\section{Introduction and Main results}
A Selberg type integral
\begin{equation}
\int_{\gamma}\prod_{1\le i<j\le m}(x_i-x_j)^g
\prod_{\substack{
1\le i\le m\\
1\le k\le n}}
(x_i-z_k)^{\lambda_{k}}
dx_1\cdots dx_m\label{(1.1)}
\end{equation}
is used to express a conformal block 
in conformal field theory 
\cite{DJMM}\cite{DF1}\cite{DF2}\cite{V} and to represent
the hypergeometric function due to Heckman and Opdam
of type $BC$ \cite{Ka}\cite{M1}. 
The integral $(1.1)$ can be thought of the pairing between 
the de Rham cohomology group
and the twisted homology group 
(the homology group with coefficients in the 
local system), which are studied by many authors, e.g. 
\cite{Ao}\cite{C}\cite{Ko}.
The study under resonance conditions, however, is
not well pursued \cite{CV}\cite{MOY}. 
The purpose of this article is to study the 
twisted homology group associated with $(\ref{(1.1)})$  
under some resonance condition.
The dimension of the space
of the regularizable cycles is determined.
\medskip

Let ${\rm Conf}_n({\Bbb C})$ denote the configuration
space of $n$ distinct points of $ {\Bbb C}:$
\smallskip

$${\rm Conf}_n({\Bbb C})=\{(z_1,\ldots, z_n)\in {\Bbb C}^n\,;\,
z_i\neq z_j\;\text{if}\;i \neq j\,\}.$$
\smallskip

\noindent
For each point $z=(z_1,\ldots, z_n)\in{\rm Conf}_n({\Bbb C}),$ 
\noindent
let $\Phi_{m,n}^{g,\lambda}$ be the function 
\begin{align*}
\Phi_{m,n}^{g,\lambda}(x;z)&=
\Phi_{m,n}^{g,\lambda_1,\ldots, \lambda_n}
(x_1,\ldots, x_m; z_1,\ldots, z_n)\\[6pt]
&=\prod_{1\le i<j\le m}(x_i-x_j)^g
\prod_{1\le i\le m}\prod_{1\le k\le n}(x_i-z_k)^{\lambda_{k}}
\end{align*} 

\noindent
on the complex manifold
\smallskip

\begin{align*}
&X_{m,n}(x_1,\ldots, x_m; z_1,\ldots, z_n)\\[6pt]
&=
\comment{{\Bbb C}^m\backslash 
\cup_{1\le i<j\le m}\;\{\,x_i-x_j=0\,\}\cup
\cup_{\substack{1\le i\le m\\1\le j\le n}}\;\{\,x_i-z_j=0\;\}.}
\left\{x\in {\Bbb C}^m\mid 
\prod_{1\le i< j\le m}(x_i-x_j)
\prod_{1\le i\le m}\prod_{1\le k\le n}(x_i-z_k)\neq 0\right\}.
\end{align*}
\smallskip

\noindent
Let ${\mathcal L}$ be the local
system determined by $\Phi_{m,n}^{g,\lambda}$.
Compactly supported twisted (loaded) chains and locally finite 
twisted (loaded) chains, together with the natural boundary operators, 
define the twisted homology group 
$H_m(X_{m,n},{\mathcal L})$ and 
the locally finite twisted homology group
$H_m^{\rm lf}(X_{m,n},{\mathcal L}),$ respectively. 
\smallskip

The symmetric group $S_m$ acts 
on $X_{m,n}$ through the coordinates $x=(x_1,\ldots, x_m).$
The action of $S_m$ on $X_{m,n}$
defines the invariant parts of $H_m(X_{m,n},{\mathcal L})$
and $H_m^{\rm lf}(X_{m,n},{\mathcal L})$, 
which will be denoted
by $H_m(X_{m,n},{\mathcal L})^{S_m}$ and 
$H_m^{\rm lf}(X_{m,n},{\mathcal L})^{S_m}.$ 
There is a natural map 
$$\iota\,:\,H_m(X_{m,n},{\mathcal L})^{S_m}\,\longrightarrow
\,H_m^{\rm  lf}(X_{m,n},{\mathcal L})^{S_m}.$$ 
A cycles in the image ${\rm Im}\,\iota$
is called a {\it regularizable cycle}, and a preimage of 
a regularizable cycle is called a regularization of it.
\medskip

If the exponent of 
the irreducible component of 
the divisor $\widetilde{D}=\pi^{-1}(D),$
where $\pi : (\widetilde{\Bbb P^1({\Bbb C})})^m\rightarrow 
({\Bbb P}^1({\Bbb C}))^m$
is the minimal blow-up along the
non-normally crossing loci of $D,$
is an ineger, the irreducible component or the 
exponent itself is said to be {\it resonant}. 
The resonance condition on the exponents makes the 
the kernel 
of the map $\iota$ non-trivial.
Hence describing the kernel and the image of the map $\iota$
under a resonance condition
is a fundamental problem.
In this paper, especially, we determine the dimension of the
space of the regularizable cycles under  
the following resonance condition:
\smallskip

\begin{equation}
2\lambda_j+g \in {\Bbb Z}\quad(1\le j\le r)
\label{eq(1.2)}
\end{equation}
\smallskip

\noindent
for each $r$ such that $0\le r\le n$. For the other divisors, 
we assume the non-resonance condition:
 
\begin{align}
&2\lambda_j+g \notin {\Bbb Z}\quad(r+1\le j\le n),
\label{eq(1.3)}\\[6pt]
&k\lambda_j+{k \choose 2}g \notin {\Bbb Z}\; 
\quad(1\le j\le n,\;k=1\;{\rm and}\;3\le k\le m)
\label{eq(1.4)}
\end{align}
\smallskip

\noindent
and 

\begin{equation}
k\lambda_\infty+{k \choose 2}g \notin {\Bbb Z}\; 
\quad(1\le k\le m),
\quad
{k \choose 2}g\notin {\Bbb Z}\quad(2\le k\le m),\label{eq(1.5)}
\end{equation}
where 

$$\lambda_\infty
=-\sum_{1\le i\le n}\lambda_i-(m-1)g
\quad{\rm and}\quad {1 \choose 2}=0\;;$$
\smallskip

\noindent
this is our assumption on 
the exponents $\lambda_j$'s and $g$, throughout this paper. 
\smallskip

In what follows, to indicate the dependence on the
number $r$, we denote
by ${\mathcal L}(r)$ the local system determined by 
$\Phi_{m,n}^{g,\lambda}$ with the resonance condition
$(\ref{eq(1.2)})$ for $r$ such that $0\le r\le n$. 
It follows from Theorem 2 of \cite{C} with
the Poincar\'e duality that, even under this resonance condition,
the rank of $H^{\rm lf}_m(X_{m,n},{\mathcal L}(r))^{S_m}$
and the rank of $H_m(X_{m,n},{\mathcal L}(r))^{S_m}$
remain to be  $D_{m,n}={n+m-2\choose m}$.
\bigskip

Let $I_{m,n}(r)$ denote the dimension of 
${\rm Im}\,\iota\subset
H_m^{\rm lf}(X_{m,n},{\mathcal L}_{m,n}(r))^{S_m},$
and $K_{m,n}(r)$ the dimension of ${\rm Ker}\,\iota\subset
H_m(X_{m,n},{\mathcal L}_{m,n}(r))^{S_m}.$ 
Then  $K_{m,n}(0)=0,\; K_{1,n}(r)=0,$ and

\begin{equation}
D_{m,n}=K_{m,n}(r)+I_{m,n}(r).
\label{eq(1.6)}
\end{equation}
\smallskip

\noindent
We have the following.
\smallskip

\noindent
{\bf Theorem 1.} {\it Suppose $(\ref{eq(1.2)}), (\ref{eq(1.3)}), 
(\ref{eq(1.4)})$ and $(\ref{eq(1.5)}).$ Then we have
\smallskip

\noindent
\begin{equation}
I_{m,n}(r)=
\left(
\begin{array}{c} 
n+m-2\\m
\end{array}\right)
{}_3F_{2}\hs
{-r,\;\dfrac{-m}{2},\;\dfrac{-m+1}{2}}
{\dfrac{-n-m+2}{2},\;\dfrac{-n-m+3}{2}}
{1},\label{eq(1.7)}
\end{equation}
\medskip

\noindent
where
\medskip

\noindent
$$
{}_3F_{2}\hs
{\alpha_1,\;\alpha_2,\;\alpha_3}
{\beta_1,\;\beta_2}
{x}
=
\sum_{k=0}^{\infty}
\frac{(\alpha_1)_k (\alpha_2)_k(\alpha_3)_k}
{(\beta_1)_k (\beta_{2})_k \,k!}\;x^k
$$
\medskip

\noindent
with $(a)_k=a(a+1)\cdots(a+k-1).$
}
\bigskip

\noindent
In the special cases $r=n$ and $r=n-1$, $(\ref{eq(1.7)})$
becomes simpler.
\medskip

\noindent
{\bf Theorem 2.} {\it Under the same condition as in Theorem 1,
we have
\smallskip

\noindent
\begin{equation*}
(1) \qquad I_{m,n}(n)={n\choose m}-{n\choose m-1},
\quad{and}\quad
(2) \qquad I_{m,n}(n-1)={n-1\choose m}.
\end{equation*}}
\smallskip

\noindent
The dimension of the
irreducible representation of the Iwahori-Hecke
algebra $H(S_n)$ parametrized by the Young diagram $(n-m,m)$
is ${n\choose m}-{n\choose m-1}$, and the representation
is realized in terms of the regularizable cycles
of $H_m^{\rm lf}(X_{m,n},{\mathcal L}_{m,n}(n))^{S_m}$ 
in \cite{M2}.
Therefore, Theorem 2  shows that the representation space 
is exactly equal to 
${\rm Im}\,\iota.$
\medskip

\section{Proof of Theorem 1}
Let $T$ be a tubular neighbourhood of 
\smallskip

\noindent
\begin{equation*}
D=\cup_{1\le i<j\le m}\{\,x_i-x_j=0 \,\}
\cup_{1\le i\le m}\cup_{1\le j\le n}
\{x_i-z_j=0\}\cup_{1\le i\le m}\{x_i=\infty\}
\end{equation*}

\noindent
in $(\PP^1(\C))^m$, and put $T^\circ=T-D.$ 
The inclusion $T^\circ\subset X_{m,n}=(\PP^1(\C))^m\backslash D$ 
leads to the exact sequence
\smallskip

\noindent
\begin{equation*}
\begin{split}
\cdots\longrightarrow 
H_{m+1}(X_{m,n},T^\circ,{\mathcal L}_{m,n}(r))^{S_m}
\longrightarrow H_m(T^\circ,{\mathcal L}_{m,n}(r))^{S_m}\\[6pt]
\longrightarrow 
H_m(X_{m,n},{\mathcal L}_{m,n}(r))^{S_m}\longrightarrow 
H_m(X_{m,n},T^\circ,{\mathcal L}_{m,n}(r))^{S_m}
\longrightarrow\cdots.
\end{split}
\end{equation*}
\smallskip

\noindent
Since the relative homology group 
$H_k(X_{m,n},T^\circ,{\mathcal L}_{m,n}(r))^{S_m}$ can be
canonically identified with 
$H_k^{\rm lf}(X_{m,n},{\mathcal L}_{m,n}(r))^{S_m}$, 
we have the exact sequence
\smallskip

\noindent
\begin{equation*}
\begin{split}
H_{m+1}^{\rm lf}(X_{m,n},{\mathcal L}_{m,n}(r))^{S_m}
\longrightarrow 
H_m(T^\circ,{\mathcal L}_{m,n}(r))^{S_m}\\[6pt]
\longrightarrow 
H_m(X_{m,n},{\mathcal L}_{m,n}(r))^{S_m}
\longrightarrow H_m^{\rm lf}(X_{m,n},{\mathcal L}_{m,n}(r))^{S_m},
\end{split}
\end{equation*}
\smallskip

\noindent
where the last arrow is the natural map $\iota.$
Theorem 2 in \cite{C} and the Poincer\'e duality imply
$H_{m+1}^{\rm lf}(X_{m,n},{\mathcal L}_{m,n}(r))^{S_m}=0.$
So our task is to study $H_m(T^\circ,{\mathcal L}_{m,n}(r))^{S_m}.$
\par\smallskip

\noindent
We first blow up $(\PP^1(\C))^m.$
Let $\pi : \widetilde{(\PP^1(\C))^m}\rightarrow (\PP^1(\C))^m$ 
be the minimal blow-up of $(\PP^1(\C))^m$ along the non-normally 
crossing loci of $D$ (cf.\,\cite{MY}).
The exponent of  

$$
\pi^{-1}\{x_{\sigma(1)}=x_{\sigma(2)}=\cdots=x_{\sigma(k)}=z_j\}
\quad(\sigma\in S_m,\; 2\le k\le m,\;1\le j\le n)$$
\smallskip

\noindent
is $k\lambda_j+{k \choose 2}g,$ that of
\smallskip

$$
\pi^{-1}\{x_{\sigma(1)}=x_{\sigma(2)}=\cdots=x_{\sigma(k)}=\infty\}
\quad(\sigma\in S_m,\; 2\le k\le m)
$$
\smallskip

\noindent
is $k\lambda_\infty+{k \choose 2}g$, and that of

\begin{equation*}
\pi^{-1}\{x_{\sigma(1)}=x_{\sigma(2)}=\cdots=x_{\sigma(k)}\} 
\quad(3\le k\le m,\;\sigma \in S_m)
\end{equation*}
\smallskip

\noindent
is ${k \choose 2}g.$
The condition $(\ref{eq(1.2)})$ means that
$\pi^{-1}\{x_{\sigma(1)}=x_{\sigma(2)}=z_j\}$ for 
$\sigma\in S_m, 1\le j\le r$ is a
resonant divisor.  \bigskip

If we denote
$\widetilde{D}=\pi^{-1}(D)$ by $\sum_j\widetilde{D}_j,$
where each $\widetilde{D}_j$ is the irreducible divisor 
(which might be an exceptional one) , a neighbourhood
of $\widetilde{D}_j$ with small radius
is set to be  $N(\widetilde{D}_j)$, and 
$$N^\circ(\widetilde{D}_j)=N(\widetilde{D}_j)
\backslash \widetilde{D}\subset\widetilde{X}_{m,n}.$$
Then we have $\widetilde{T}^\circ=\sum_jN^\circ(\widetilde{D}_j)
\subset\widetilde{X}_{m,n}.$
\medskip

If $\widetilde{D}_i$ is an irreducible non-resonant 
divisor, we have

$$H_k(N^\circ(\widetilde{D}_i),\pi^*{\mathcal L}_{m,n}(r) )=0
\qquad (0\le k\le 2m).$$
\smallskip

\noindent
If $\widetilde{D}_i$ is an irreducible non-resonant divisor,
and if $\widetilde{D}_j$ is an irreducible resonant divisor, 
we have

$$H_k(N^\circ(\widetilde{D}_i)\cap N^\circ(\widetilde{D}_j),
\pi^*{\mathcal L}_{m,n}(r) )=0\qquad (0\le k\le 2m).$$
\smallskip

\noindent
This is because $N^\circ(\widetilde{D}_i)$ and 
$N^\circ(\widetilde{D}_i)\cap N^\circ(\widetilde{D}_j)$
have a punctured $1$-disc, as a direct summand, carrying
a non-resonant exponent.
\medskip

Hence,  
the Mayer-Vietoris sequence 

\begin{equation*}
\begin{split}
&\cdots \longrightarrow
H_m(U\cap V,\pi^*{\mathcal L}_{m,n}(r))\longrightarrow 
H_m(U,\pi^*{\mathcal L}_{m,n}(r))\oplus 
H_m(V,\pi^*{\mathcal L}_{m,n}(r))\\[6pt]
&\longrightarrow H_m(U\cup V,\pi^*{\mathcal L}_{m,n}(r))
\longrightarrow H_{m-1}(U\cap V,\pi^*{\mathcal L}_{m,n}(r))
\longrightarrow\cdots,
\end{split}
\end{equation*}
\smallskip

\noindent
for $U=N^\circ(\widetilde{D}_i)$ and 
$V=N^\circ(\widetilde{D}_j)$ stated above,
implies 

$$H_m(U\cup V,\pi^*{\mathcal L}_{m,n}(r))
\cong
H_m(V,\pi^*{\mathcal L}_{m,n}(r)).
$$
\smallskip

\noindent
Applying this equivalence repeatedly, 
we have 

\begin{align}
&H_m(T^\circ,{\mathcal L}_{m,n}(r))
\cong H_m(\widetilde{T}^\circ,\pi^*{\mathcal L}_{m,n}(r))
\nonumber\\[5pt]
\cong&H_m(\cup_{1\le k\le r}\cup_{1\le i<j\le m}
N^\circ(\pi^{-1}\{x_i=x_j=z_k\}),
\pi^*{\mathcal L}_{m,n}(r)).
\end{align}
\smallskip

\noindent
Thus

\begin{align}
&H_m(T^\circ,{\mathcal L}_{m,n}(r))^{S_m}
\cong H_m(\widetilde{T}^\circ,\pi^*{\mathcal L}_{m,n}(r))^{S_m}
\nonumber\\[5pt]
\cong&H_m(\cup_{1\le k\le r}\cup_{1\le i<j\le m}
N^\circ(\pi^{-1}\{x_i=x_j=z_k\}),
\pi^*{\mathcal L}_{m,n}(r))^{S_m}.
\end{align}

\medskip

\noindent
At this stage, we have the following.
\begin{lem}
For $m\ge 2,$
\begin{equation*}
K_{m,n}(r)=D_{m-2,n}+K_{m,n}(r-1)-K_{m-2,n}(r-1),
\end{equation*}
where $D_{0,n}=0.$
\end{lem}
\noindent
{\rm\bf Proof.}
Set 

$$U=\cup_{k=1}^{r-1}\cup_{\sigma\in S_m}
N^\circ(\pi^{-1}
\{x_{\sigma(1)}=x_{\sigma(2)}=z_k\})$$
\smallskip

\noindent
and

$$V=\cup_{\sigma\in S_m}
N^\circ(\pi^{-1}
\{x_{\sigma(1)}=x_{\sigma(2)}=z_r\}).$$
\smallskip

\noindent
Then 

\begin{equation}
\begin{split}
&0\longrightarrow 
H_m(U\cap V,\pi^*{\mathcal L}_{m,n}(r))^{S_m}\longrightarrow 
H_m(U,\pi^*{\mathcal L}_{m,n}(r))^{S_m}
\oplus H_m(V,\pi^*{\mathcal L}_{m,n}(r))^{S_m}\\[6pt]
&\longrightarrow H_m(U\cup V,\pi^*{\mathcal L}_{m,n}(r))^{S_m}
\longrightarrow 0,
\end{split}
\end{equation}
\smallskip

\noindent
since

\begin{equation*}
H_{k}(U\cap V,\pi^*{\mathcal L}_{m,n}(r))=0
\qquad\text{for}\qquad k\neq m.
\end{equation*}
\medskip

By definition, we have
 
\begin{equation}
{\rm dim}\,H_m(U\cup V,\pi^*{\mathcal L}_{m,n}(r))^{S_m}
=K_{m,n}(r),
\end{equation}
and  
\begin{equation}
{\rm dim}\,H_m(U,\pi^*{\mathcal L}_{m,n}(r))^{S_m}
=K_{m,n}(r-1).
\end{equation}
\smallskip

\noindent
Lemma 3 stated in Section 4 reads

\begin{equation*}
H_m(V,\pi^*{\mathcal L}_{m,n}(r))^{S_m}
\cong
H_{m-2}(X_{m-2,n},{\mathcal L}'_{m-2,n}(r-1))^{S_{m-2}},
\end{equation*}
\smallskip

\noindent
where ${\mathcal L}'_{m-2,n}(r-1)={\mathcal L}_{m-2,n}
|_{\lambda_r\to\lambda_r+2g}(r-1)$ is the local system 
on $X_{m-2,n}$ determined by

$$\Phi_{m-2,n}^{g,\lambda}(x;z)\prod_{i=1}^{m-2}(x_i-z_r)^{2g}.$$
\smallskip

\noindent
Hence, we have

\begin{equation}
{\rm dim}\,H_m(V,\pi^*{\mathcal L}_{m,n}(r))^{S_m}=D_{m-2,n}.
\end{equation}
\smallskip

\noindent
On the other hand, we have

\begin{align*}
&U\cap V\\
\comment{&=\cup_{k=1}^{r-1}\cup_{\sigma\in S_m}\biggl\{
N^\circ\Bigl(\pi^{-1}
\{x_{\sigma(1)}=x_{\sigma(2)}=z_k\}\Bigr)
\cap
N^\circ\Bigl(\pi^{-1}
\{x_{\sigma(3)}=x_{\sigma(4)}=z_r\}\Bigr)\biggr\}\\}
&=\cup_{\sigma\in S_m}\biggl\{\left(
\cup_{k=1}^{r-1}\;
N^\circ\Bigl(\pi^{-1}
\{x_{\sigma(1)}=x_{\sigma(2)}=z_k\}\Bigr)\right)
\cap
N^\circ\Bigl(\pi^{-1}
\{x_{\sigma(3)}=x_{\sigma(4)}=z_r\}\Bigr)\biggr\}.
\end{align*}
\smallskip

\noindent
Hence the same argument for proving Lemma 3 shows

\begin{align*}
&H_m(U\cap V,\pi^*{\mathcal L}_{m,n}(r))^{S_m}\\[6pt]
&\cong 
H_{m-2}(
\cup_{k=1}^{r-1}\cup_{\sigma\in S_{m-2}}
N^\circ(\pi^{-1}
\{x_{\sigma(1)}=x_{\sigma(2)}=z_k\})
,\pi^*{\mathcal L}'_{m-2,n}(r-1))^{S_{m-2}},
\end{align*}
\smallskip

\noindent
where $\pi$ in the righthand-side is the restriction of
$\pi : \widetilde{(\PP^1(\C))^m}\rightarrow (\PP^1(\C))^m$ 
on the $m-2$ space: $x_{m-1}=x_m=0.$
\medskip

\noindent
Therefore, we have

\begin{equation}
{\rm dim}\,H_m(U\cap V,\pi^*{\mathcal L}_{m,n}(r))^{S_m}
=K_{m-2,n}(r-1).
\end{equation}
\smallskip

\noindent
Combination of $(2.2)$ and $(2.3-6)$ implies the required result.
\begin{flushright}
$\square$
\end{flushright}
\medskip

Furthermore we have

\begin{prop}
For $m\ge 2,$

\begin{equation*}
K_{m,n}(r)=rD_{m-2,n}-K_{m-2,n}(1)-K_{m-2,n}(2)-
\cdots-K_{m-2,n}(r-1).
\end{equation*}
\smallskip

\end{prop}
\noindent
{\rm\bf Proof.}
Summing up the equalities

\begin{align*}
K_{m,n}(r)&=D_{m-2,n}+K_{m,n}(r-1)-K_{m-2,n}(r-1),\\[5pt]
K_{m,n}(r-1)&=D_{m-2,n}+K_{m,n}(r-2)-K_{m-2,n}(r-2),\\[5pt]
&\cdots\\
K_{m,n}(3)&=D_{m-2,n}+K_{m,n}(2)-K_{m-2,n}(2),\\[5pt]
K_{m,n}(2)&=D_{m-2,n}+K_{m,n}(1)-K_{m-2,n}(1),
\end{align*}
\smallskip

\noindent
each of which is a special case of the equality in Lemma 1,
we have the required result, since $K_{m,n}(1)=K_{m-2,n}(1)=0.$
\begin{flushright}
$\square$
\end{flushright}
\smallskip

Theorem 1.3 of \cite{MOY} implies $K_{2,n}(r)=r.$ 
Proposition 1 implies

\begin{equation}
K_{3,n}(r)=rD_{1,n},
\end{equation}
\smallskip

\noindent
since $K_{1,n}(1)=\cdots=K_{1,n}(r-1)=0.$
Generally we have the following.

\begin{prop}For $m\ge 2$ and $r\ge 0$,

\begin{align}
K_{m,n}(r)&={r \choose 1}D_{m-2,n}-
{r \choose 2}D_{m-4,n}+\cdots
+(-1)^{s-1}{r \choose s}D_{m-2s,n}+\cdots\nonumber\\[5pt]
&=\sum_{s\ge 1}(-1)^{s-1}{r \choose s}D_{m-2s,n},
\end{align}

\noindent
where ${r \choose s}=0$ for $s>r.$
\end{prop}
\noindent
{\bf\rm Proof.}
Suppose the equality in case $m.$ Then 
Proposition 1 shows

\begin{align*}
K_{m+2,n}(r)&=rD_{m,n}-K_{m,n}(1)-K_{m,n}(2)-
\cdots-K_{m,n}(r-1)\\[5pt]
&=rD_{m,n}-\sum_{1\le t\le r-1}\sum_{1\le s}
(-1)^{s-1}{t \choose s}D_{m-2s,n}\\[5pt]
&=rD_{m,n}-\sum_{1\le s}(-1)^{s-1}
\sum_{1\le t\le r-1}
{t \choose s}D_{m-2s,n}\\[5pt]
&=rD_{m,n}-\sum_{1\le s}(-1)^{s-1}
{r \choose s+1}D_{m-2s,n}\\
&=\sum_{s\ge 1}(-1)^{s-1}{r \choose s}D_{m+2-2s,n}.
\end{align*}
\smallskip

\noindent
For the fourth equality, we used Lemma 2 below. 
The induction on $m$ leads to the required result. 
$\square$

\begin{lem}
$$1+{s+1\choose s}+{s+2\choose s}+\cdots+{r-1\choose s}
={r\choose s+1}.$$
\end{lem}

\noindent
{\rm\bf Proof.}

\begin{align*}
{r\choose s+1}&={r-1\choose s+1}+{r-1\choose s}\\[5pt]
&=\left\{{r-2\choose s+1}+{r-2\choose s}\right\}+{r-1\choose s}\\[5pt]
&=\left\{{r-3\choose s+1}+{r-3\choose s}\right\}
+{r-2\choose s}+{r-1\choose s}\\[5pt]
&\cdots\\
&=\left\{{s+1\choose s+1}+{s+1\choose s}\right\}
+{s+2\choose s}+\cdots+{r-1\choose s}\\[5pt]
&=1+{s+1\choose s}+{s+2\choose s}+\cdots
+{r-1\choose s}. 
\end{align*}

\begin{flushright}
$\square$
\end{flushright}

\bigskip

By Proposition 2 and $(1.6)$, we obtain 

\begin{align}
&I_{m,n}(r)=D_{m,n}-K_{m,n}(r)\nonumber\\[5pt]
&=\sum_{s=0}^{[m/2]}(-1)^s{r \choose s}
D_{m-2s,n},\label{eq(2.10)}
\end{align}

\noindent
where $[x]$ denotes the largest integer not exceeding $x$. 
\bigskip

Finally, we express $(2.10)$ 
in terms of the generalized hypergeometric series.
The equalities

\begin{equation*}
(n+m-2)!
=2^{2s}((-n-m+2)/2)_s
((-n-m+3)/{2})_s(n+m-2-2s)!
\end{equation*}

\noindent
and 

$$m!=2^{2s}(-m/2)_s
((-m+1)/2)_s(m-2s)!$$
\smallskip

\noindent
lead to

$${n+m-2-2s\choose m-2s}={n+m-2\choose m}
\frac{(-m/2)_s((-m+1)/2)_s}
{((-n-m+2)/2)_s((-n-m+3)/2)_s}.
$$
\smallskip

\noindent
On the other hand, 
$${r\choose s}=(-1)^s
\frac{(-r)_s}{s!}.
$$
\medskip

\noindent
Hence, we reach the desired expression
\smallskip

$$I_{m,n}(r)=
\left(
\begin{array}{c} 
n+m-2\\m
\end{array}\right)
{}_3F_{2}\hs
{-r,\;-m/2,\;(-m+1)/2}
{(-n-m+2)/2,\;(-n-m+3)/2}
{1}.
$$
\bigskip

\noindent
This completes the proof of Theorem 1.

\section{Proof of Theorem 2}
We first derive $(1)$ of Theorem 2. 
Note that,
in case $m$ is even $2j$, 

\begin{equation*}
I_{m,n}(n)=
\left(
\begin{array}{c} 
n+2j-2\\2j
\end{array}\right)
{}_3F_{2}\hs
{-n,\;-j,\;\dfrac{1}{2}-j}
{-j+1-\dfrac{n}{2},\;-j+\dfrac{3-n}{2}}
{1},
\end{equation*}
\smallskip

\noindent
and, in case $m$ is odd $2j+1$, 

\begin{equation*}
I_{m,n}(n)=
\left(
\begin{array}{c} 
n+2j-1\\2j+1
\end{array}\right)
{}_3F_{2}\hs
{-n,\;-j-\dfrac{1}{2},\;-j}
{\dfrac{-n+1}{2}-j,\;1-j-\dfrac{n}{2}}
{1}.
\end{equation*}
\medskip

Pfaff-Saalsch\"utz's theorem (Theorem 2.2.6 in \cite{AAR})

\begin{equation}
{}_3F_{2}\hs
{\;a,\;b,\;-j}
{\;c,\;1+a+b-c-j\;}
{1}
=\frac{(c-a)_j(c-b)_j}{(c)_j(c-a-b)_j}\qquad
{\rm for}\quad j\in {\Bbb Z}_{\ge 0}
\end{equation}
\smallskip

\noindent
and a contiguity relation

\begin{align*}
&(b-a)\,{}_3F_{2}\hs
{\;a,\;b,\;-j\;}
{\;c,\;a+b-c+2-j\;}
{1}+
a\,{}_3F_{2}\hs
{\;a+1,\;b,\;-j\;}
{\;c,\;a+b-c+2-j\;}
{1}\nonumber\\[6pt]
&-b\,{}_3F_{2}\hs
{\;a,\;b+1,\;-j\;}
{\;c,\;a+b-c+2-j\;}
{1}=0,  
\end{align*}
\smallskip

\noindent
which follows from the identity

$$
a(a+1)_k(b)_k-b(a)_k(b+1)_k=(a-b)(a)_k(b)_k,
$$
\smallskip

\noindent
imply
\begin{align*}
&{}_3F_{2}\hs
{\;a,\;b,\;-j\;}
{\;c,\;a+b-c+2-j\;}
{1}
\nonumber\\[6pt]
&=\frac{a}{a-b}
\frac{(c-a-1)_j(c-b)_j}{(c)_j(c-a-b-1)_j}
+\frac{b}{b-a}
\frac{(c-a)_j(c-b-1)_j}{(c)_j(c-a-b-1)_j}.
\end{align*}
\medskip

\noindent
Hence we have
\begin{align*}
&{}_3F_{2}\hs
{-n,\;-j,\;\dfrac{1}{2}-j}
{-j+1-\dfrac{n}{2},\;-j+\dfrac{3-n}{2}}
{1}\\[8pt]
&=\frac{-n}{-n+j-\dfrac{1}{2}}\;
\frac{(1-\dfrac{n}{2})_j(\dfrac{1-n}{2})_j}
{(\dfrac{n}{2})_j(\dfrac{n-1}{2})_j}
+
\frac{-j+\dfrac{1}{2}}{n-j+\dfrac{1}{2}}\;
\frac{(-\dfrac{n}{2})_j(-\dfrac{1+n}{2})_j}
{(\dfrac{n}{2})_j(\dfrac{n-1}{2})_j},
\end{align*}
\smallskip

\noindent
which turns out to be

\begin{align*}
&\frac{1}{2}(-1+4j-n)
\frac{(-\dfrac{n}{2})_j(\dfrac{1-n}{2})_{j-1}}
{(\dfrac{n}{2})_j(\dfrac{n-1}{2})_j}\\[8pt]
&=(n+1-4j)
\frac{n(n-1)\cdots(n-2j+2)}
{(n+2j-2)(n+2j-3)\cdots n(n-1)}.
\end{align*}
\smallskip

\noindent
Therefore, we have

$$I_{2j,n}(n)=\frac{n(n-1)\cdots(n-2j+2)}{(2j)!}(n+1-4j).$$
\medskip

In the same way, we have

\begin{align*}
&{}_3F_{2}\hs
{-n,\;-j-\dfrac{1}{2},\;-j}
{\dfrac{-n+1}{2}-j,\;1-j-\dfrac{n}{2}}
{1}\\[8pt]
&=(n+1-4j)
\frac{(n-2)(n-3)\cdots(n-2j+1)}
{(n+2j-1)(n+2j-3)\cdots (n+1)},
\end{align*}
\smallskip

\noindent
and thus

$$
I_{2j+1,n}(n)=\frac{n(n-1)\cdots(n-2j+1)}{(2j+1)!}(n-4j-1).
$$
\medskip

\noindent
This completes the proof of (1).
\bigskip

Next we derive (2) of Theorem 2.
Note that

\begin{equation*}
I_{2j,n}(n-1)=
\left(
\begin{array}{c} 
n+2j-2\\2j
\end{array}\right)
{}_3F_{2}\hs
{-n+1,\;-j,\;\dfrac{1}{2}-j}
{-j+1-\dfrac{n}{2},\;-j+\dfrac{3-n}{2}}
{1}
\end{equation*}
\smallskip

\noindent
and

\begin{equation*}
I_{2j+1,n}(n-1)=
\left(
\begin{array}{c} 
n+2j-1\\2j+1
\end{array}\right)
{}_3F_{2}\hs
{-n+1,\;-j-\dfrac{1}{2},\;-j}
{\dfrac{-n+1}{2}-j,\;1-j-\dfrac{n}{2}}
{1}.
\end{equation*}
\smallskip

\noindent
Pfaff-Saalsch\"utz's theorem $(3.1)$ implies

\begin{align*}
&{}_3F_{2}\hs
{-n+1,\;-j,\;\dfrac{1}{2}-j}
{-j+1-\dfrac{n}{2},\;-j+\dfrac{3-n}{2}}
{1}
=\frac{(\dfrac{n}{2}-j)_j(\dfrac{1-n}{2})_j}
{(1-\dfrac{n}{2}-j)_j(\dfrac{1-n}{2})_j}\\[8pt]
&=\frac{(1-\dfrac{n}{2})_j(\dfrac{1-n}{2})_j}
{(\dfrac{n}{2})_j(\dfrac{n-1}{2})_j}
=\frac{(1-n)_{2j}}{(n-1)_{2j}}
\end{align*}
\smallskip

\noindent
and

\begin{equation*}
{}_3F_{2}\hs
{-n+1,\;-j-\dfrac{1}{2},\;-j}
{\dfrac{-n+1}{2}-j,\;1-j-\dfrac{n}{2}}
{1}
=\frac{(\dfrac{3-n}{2})_j(1-\dfrac{n}{2})_j}
{(\dfrac{n+1}{2})_j(\dfrac{n}{2})_j}
=-\frac{(1-n)_{2j+1}}{(n-1)_{2j+1}}.
\end{equation*}
\smallskip

\noindent
Therefore, we obtain

$$I_{2j,n}(n-1)=\frac{(1-n)_{2j}}{(1)_{2j}}={n-1\choose 2j}$$
\smallskip

\noindent
and

$$I_{2j+1,n}(n-1)=\frac{(1-n)_{2j+1}}{(1)_{2j+1}}
={n-1\choose 2j+1}.$$
\medskip

\noindent
This completes the proof of $(2)$.

\section{A Key Lemma}

In this section, we fix $z_1$ to be $0$ for simplicity, and
we study the homological structure around the
subvariety $\{(x_1,\ldots,x_m)\in {\Bbb C}^m\mid
x_i=x_j=0\}$ of codimension two for  $1\le i<j\le m$. 
\smallskip

Choose an appropriate open set $N^\circ(x_i=x_j=0)\subset
X_{m,n}$ such that 
$H_m(N^\circ(x_i=x_j=0), {\mathcal L}_{m,n}(r)) \cong
H_m(N^\circ(\pi^{-1}\{x_i=x_j=0\}), 
\pi^*{\mathcal L}_{m,n}(r))$. 
For example,  we take

\begin{equation*}
\begin{split}
&N^\circ(x_i=x_j=0)\\[6pt]
&:=\{(x_1,\ldots,x_m)
\in X_{m,n} \mid |x_i|,|x_j|<\e(x_1,\ldots,\hat{x_i},
\ldots,\hat{x_j},\ldots,x_m)\},
\end{split}
\end{equation*}
\smallskip

\noindent
where 

\begin{equation*}
\e(x_3,\ldots,x_m)
=\min\{\e|x_3|,\ldots, \e|x_m|, \e^2\},
\end{equation*}
\smallskip

\noindent
and $\e$ is a positive number satisfying 
$\e<\min\{1,|z_2|,\ldots, |z_n|\}.$

\begin{lem} Suppose $(\ref{eq(1.2)}).$
Then
we have the isomorphism

\begin{equation*}
\begin{split}
&H_m(\cup_{1\le i<j\le m}
N^\circ(x_{i}=x_{j}=0),
{\mathcal L}_{m,n}(r))^{S_m}\\[6pt]
&\cong 
H_{m-2}(X_{m-2,n},{\mathcal L}_{m-2,n}'(r-1))^{S_{m-2}}, 
\end{split}
\end{equation*}
\smallskip

\noindent
where the local system ${\mathcal L}_{m-2,n}'(r-1)$
is determined by

$$
\prod_{1\le i<j\le m-2}(x_i-x_j)^g
\prod_{1\le i\le m-2}
\left\{x_i^{\lambda_{1}+2g}
\prod_{2\le k\le n}(x_i-z_k)^{\lambda_{k}}
\right\}.$$ 
\end{lem}
\medskip

\noindent
{\bf Proof.} 
The function $\Phi_{m,n}^{g,\lambda}$ defining the local system 
${\mathcal L}_{m,n}(r)$ 
restricted on $N^\circ(x_1=x_2=0)$ is reduced to 

\begin{equation}
\prod_{1\le i<j\le m}(x_i-x_j)^g\prod_{1\le i\le m}
x_i^{\lambda_1}
\prod_{3\le i\le m}\prod_{2\le k\le n}(x_i-z_k)^{\lambda_k}.
\end{equation}
\smallskip

\noindent
The pullback of $(4.1)$ by the morphism 

$$
\psi:(s_1,s_2,x_3,\ldots,x_m)\to(x_1,\ldots,x_m),
$$

\noindent
where

$$\quad x_1=s_1x_3\cdots x_m, \quad 
x_2=s_2x_3\cdots x_m$$
\smallskip

\noindent
is expressed by

\begin{align*}
&(s_1s_2)^{\lambda_1}
(s_1-s_2)^g\\[6pt]
\times 
&\prod_{3\le i<j\le m}(x_i-x_j)^g
\prod_{3\le i\le m}
\left\{x_i^{3\lambda_1+3g}
\left(1-s_1\frac{x_3\cdots x_m}{x_i}\right)^g
\left(1-s_2\frac{x_3\cdots x_m}{x_i}\right)^g
\right\}\\[6pt]
\times
&\prod_{3\le i\le m}\prod_{2\le k\le n}(x_i-z_k)^{\lambda_k}.
\end{align*}
\smallskip

On the other hand, the condition
$|x_1|<\e|x_i|$ and $|x_2|<\e|x_i|$
on $N^\circ(x_1=x_2=0)$ implies 
$$\left|s_1\dfrac{x_3\cdots x_m}{x_i}\right|<\e
\quad\text{and}\quad
\left|s_2\dfrac{x_3\cdots x_m}{x_i}\right|<\e$$
on $\psi^{-1}N^\circ(x_1=x_2=0),$
for $3\le i\le m$. 
This shows that  
the contribution of the factors

$$\prod_{3\le i\le m}
\left(1-s_1\frac{x_3\cdots x_m}{x_i}\right)^g
\left(1-s_2\frac{x_3\cdots x_m}{x_i}\right)^g
$$
\smallskip

\noindent
to the local sytem $\psi^*{\mathcal L}_{m,n}(r)$ is trivial.
\medskip

Thus, the local system  $\psi^*{\mathcal L}_{m,n}(r)$ on 
$\psi^{-1}N^\circ(x_1=x_2=0)$
can be considered as that determined by
\medskip

$$
(s_1s_2)^{\lambda_1}
(s_1-s_2)^g
\times 
\prod_{3\le i<j\le m}(x_i-x_j)^g
\prod_{3\le i\le m}
x_i^{\lambda_1+2g}
\prod_{3\le i\le m}\prod_{2\le k\le n}(x_i-z_k)^{\lambda_k},
$$
\medskip

\noindent
(Note that  $2\lambda_1+g\in {\Bbb Z}$). 
This fact leads to the isomorphism
\medskip

\noindent
\begin{equation*}
\begin{split}
&H_m(N^\circ(x_1=x_2=0),{\mathcal L}_{m,n}(r))\\[6pt]
&\cong 
H_2(X_{2,1}(x_1,x_2;0),{\mathcal L}_{2,1}(1))\otimes 
H_{m-2}(X_{m-2,n}(x_3, \ldots, x_m;z),{\mathcal L}'_{m-2,n}(r-1)).
\end{split}
\end{equation*}

Since ${\rm rank}\,H_2(X_{2,1},{\mathcal L}_{2,1}(1))=1$
from Theorem 1.3 of \cite{MOY}, we fix
a generator of $H_2(X_{2,1},{\mathcal L}_{2,1}(1))$
and identify $H_2(X_{2,1},{\mathcal L}_{2,1}(1))$ 
with ${\Bbb C}.$ Then we have

\begin{equation*}
\begin{split}
&H_m(N^\circ(x_i=x_j=0),{\mathcal L}_{m,n}(r))\\[6pt]
&\cong 
H_2(X_{2,1}(x_i,x_j;0),{\mathcal L}_{2,1}(1))\\[6pt]
&\qquad\otimes 
H_{m-2}(X_{m-2,n}
(x_1, \ldots, \hat{x_i}, \ldots,  \hat{x_j}, \ldots, x_m;z),
{\mathcal L}'_{m-2,n}(r-1))\\[6pt]
&\cong 
H_{m-2}(X_{m-2,n}
(x_1, \ldots, \hat{x_i}, \ldots,  \hat{x_j}, \ldots, x_m;z),
{\mathcal L}'_{m-2,n}(r-1))\\[6pt]
&\cong 
H_{m-2}(X_{m-2,n}
(x_1, \ldots, x_{m-2};z),
{\mathcal L}'_{m-2,n}(r-1))
\end{split}
\end{equation*}
\medskip

\noindent
for $1\le i<j\le m$. Henceforce,

\begin{equation*}
\begin{split}
&H_m(\cup_{1\le i<j\le m}
N^\circ(x_{i}=x_{j}=0),
{\mathcal L}_{m,n}(r))^{S_m}\\[6pt]
&\cong 
H_{m-2}(X_{m-2,n},{\mathcal L}'(r-1))^{S_{m-2}}. 
\end{split}
\end{equation*}
\medskip

\noindent
This completes the proof.\qquad$\square$

\bigskip
\begin{flushleft}
Katsuhisa Mimachi\\
Department of Mathematics\\
Tokyo Institute of Technology\\
Oh-okayama, Meguro-ku, Tokyo 152-8551\\
Japan\\
mimachi@math.titech.ac.jp
\bigskip\\

Masaaki Yoshida\\
Department of Mathematics\\
Kyushu University\\
Ropponmatsu, Fukuoka 810-8560\\
Japan\\
myoshida@math.kyushu-u.ac.jp
\bigskip
\end{flushleft}

\end{document}